# CONTINUUM PERCOLATION WITH STEPS IN AN ANNULUS

By Paul Balister[1], Béla Bollobás[2] and Mark Walters[1]

*University of Memphis*

Let $A$ be the annulus in $\mathbb{R}^2$ centered at the origin with inner and outer radii $r(1-\varepsilon)$ and $r$, respectively. Place points $\{x_i\}$ in $\mathbb{R}^2$ according to a Poisson process with intensity 1 and let $\mathcal{G}_A$ be the random graph with vertex set $\{x_i\}$ and edges $x_i x_j$ whenever $x_i - x_j \in A$. We show that if the area of $A$ is large, then $\mathcal{G}_A$ almost surely has an infinite component. Moreover, if we fix $\varepsilon$, increase $r$ and let $n_c = n_c(\varepsilon)$ be the area of $A$ when this infinite component appears, then $n_c \to 1$ as $\varepsilon \to 0$. This is in contrast to the case of a "square" annulus where we show that $n_c$ is bounded away from 1.

**1. Introduction.** Consider the following percolation process. Fix $r > 0$ and $\varepsilon$ with $0 \leq \varepsilon \leq 1$. Place points $\{x_i\}$ in $\mathbb{R}^2$ according to a Poisson process with intensity 1. Join points $x_i$ and $x_j$ provided the Euclidean distance $\|x_i - x_j\|_2$ lie between $r(1-\varepsilon)$ and $r$. We wish to know whether the resulting graph has an infinite component.

More, generally, let $A$ be a bounded region in $\mathbb{R}^2$ and define an infinite random graph $\mathcal{G}_A$ with vertex set $V(\mathcal{G}_A) \subseteq \mathbb{R}^2$ given by a Poisson process with intensity 1 and edges $xy \in E(\mathcal{G}_A)$ when $x \in y + A$ (where $y + A = \{y + a : a \in A\}$). Define, for $x \in \mathbb{R}^2$, $A_\varepsilon(x, r) = \{y : (1-\varepsilon)r \leq \|y - x\|_2 \leq r\}$ to be the annulus centered at $x$ with inner radius $r(1-\varepsilon)$ and outer radius $r$. Then the percolation process above is given by $\mathcal{G}_A$, where $A = A_\varepsilon(0, r)$.

For convex $A$, the model $\mathcal{G}_A$ has been widely studied; see, for example, Grimmett (1999), Häggström and Meester (1996), Jonasson (2001) and especially Meester and Roy (1996). However, very little is known about the model for non-convex bodies such as the annulus $A_\varepsilon$.

The models $\mathcal{G}_A$ occur very naturally as ad hoc wireless networks. Suppose that transmitters are distributed at random in the plane and broadcast to

Received June 2003; revised September 2003.
[1]Supported by NSF Grant EIA-01-30352.
[2]Supported by NSF Grants DMS-99-70404 and EIA-01-30352 and DARPA Grant F33615-01-C1900.

neighbors if they are inside a certain region relative to the transmitter (i.e., the set $A$). It is then very natural to ask whether a message can propogate through the network (i.e., does an infinite connected component exist?).

The model $\mathcal{G}_{A_\varepsilon(r)}$ is monotonic in $r$ since, by scaling, it is equivalent to one in which $r$ (and hence $A$) is fixed and the intensity of the Poisson process varies. Hence, by standard results, for any $\varepsilon$ there is a critical value $n_c(\varepsilon)$ of the area $|A| = \pi r^2 \varepsilon (2-\varepsilon)$ (or, equivalently, the radius $r$) above which an infinite component exists with probability 1 and below which there is almost surely no infinite component.

Our main result is the following:

THEOREM 1.
$$1 + \frac{\varepsilon}{\pi\sqrt{3}} \leq n_c(\varepsilon) \leq 1 + o(\varepsilon^{1/15}).$$

In particular, the critical area tends to 1 as $\varepsilon \to 0$. The area $|A|$ is just the expected number of neighbors of a point, so for small $\varepsilon$, an infinite component appears when the average degree of $\mathcal{G}_A$ is slightly more than 1. This is not quite as surprising as it may seem: The giant component in a random graph appears when the average degree of the graph is 1 [see, e.g., Bollobás (2001)].

By comparison, if instead of the annulus $A = A_\varepsilon(0, r)$ we use the "square" annulus defined by the $l_\infty$-norm, then, as we show, $n_c(\varepsilon)$ is bounded away from 1 independently of $\varepsilon$.

We abuse notation by writing $|S|$ for the standard Lebesgue measure in either $\mathbb{R}^2$ or $\mathbb{R}$, or for the number of elements of $S$, depending on context.

The layout of this article is as follows: In the next section we prove the lower bound for Theorem 1 and the lower bound for the square annulus. In Section 3 we prove the upper bound in Theorem 1 modulo a key proposition (Proposition 6). In the final section we prove Proposition 6.

After this work was done, we found out that, independently and simultaneously, Franceschetti, Booth, Cook, Meester and Bruck (2003) proved the same result.

**2. Lower bound.** We first prove the lower bound on $n_c(\varepsilon)$ from Theorem 1. We start with a simple geometric lemma.

LEMMA 2. *Suppose* $y - x \in A = A_\varepsilon(0, r)$. *Then* $|A_\varepsilon(x, r) \cap A_\varepsilon(y, r)| \geq \frac{\varepsilon}{\pi\sqrt{3}}|A|$.

PROOF. The proof is easy but tedious, so we shall not spell it out. For small $\varepsilon$, the worst case occurs when $\|x - y\|_2 = r$. □

PROOF OF THE LOWER BOUND OF THEOREM 1. Suppose $|A| \leq 1 + \frac{\varepsilon}{\pi\sqrt{3}}$ and fix a vertex $x_1$ of the Poisson process. Let $E_n$ be the expected number of induced paths of length $n$ starting at $x_1$. We show that $E_n \to 0$ as $n \to \infty$. Suppose $n \geq 2$ and we have an induced path $x_1, \ldots, x_n$ of length $n$. If this is the initial segment of an induced path $x_1, \ldots, x_{n+1}$ of length $n+1$, then $x_{n+1} \in A_\varepsilon(x_n, r) \setminus \bigcup_{i=1}^{n-1} A_\varepsilon(x_i, r)$, so the expected number of points $x_{n+1}$, conditioned on the points $x_1, \ldots, x_n$ and on the points in $\bigcup_{i=1}^{n-1} A_\varepsilon(x_i, r)$ is at most $|A_\varepsilon(x_n, r) \setminus \bigcup_{i=1}^{n-1} A_\varepsilon(x_i, r)| \leq |A_\varepsilon(x_n, r) \setminus A_\varepsilon(x_{n-1}, r)|$. However, $x_n - x_{n-1} \in A$, so by Lemma 2, this is at most $\alpha = |A|(1 - \frac{\varepsilon}{\pi\sqrt{3}}) < 1$. Hence the expected number of points $x_{n+1}$ conditioned on the points $x_1, \ldots, x_n$ is at most $\alpha$ and $E_{n+1} \leq \alpha E_n$. It follows that $E_n \to 0$ as $n \to \infty$, but

$$\mathbb{P}(x_1 \text{ is in an infinite component}) \leq \mathbb{P}(\exists \text{ induced path } x_1, \ldots, x_n) \leq E_n$$

and hence, with probability 1, $x_1$ is not in an infinite component. Since $x_1$ was chosen arbitrarily and there are almost surely countably many points in the Poisson process, $\mathcal{G}_A$ almost surely has no infinite component. □

We show in the next section that $n_c \to 1$ as $\varepsilon \to 0$. However, first we show that this is not true for the square annulus, defined by the $l_\infty$-norm.

LEMMA 3. *Let $A$ be a region in $\mathbb{R}^2$ with $A = -A$. Define a percolation process by joining $x_i$ and $x_j$ when $x_i - x_j \in A$. Then the Poisson process with density 1 almost surely does not have an infinite component if $|A|^3 - \int_{A \times A} |(x + A) \cap (y + A)| \, dx \, dy < 1$.*

PROOF. As before, we count the expected number of induced paths of length $n$ starting at $x_1$. Suppose $x_1, \ldots, x_n$ are fixed. We count the expected number of choices of $x_{n+1}, x_{n+2}$ and $x_{n+3}$ that give an induced path $x_1, \ldots, x_{n+3}$, conditional on $x_1, \ldots, x_n$. Write $x_{n+1} = x_n + x$ and $x_{n+2} = x_{n+1} - y$, so that $x, y \in A = -A$. Now $x_{n+3} \in A' = (x_{n+2} + A) \setminus (x_n + A)$. The area of this set is $|A| - |(x + A) \cap (y + A)|$. From this we see that

$$E_{n+3} \leq \left( \int_{A \times A} |A'| \, dx \, dy \right) E_n = \left( |A|^3 - \int_{A \times A} |(x + A) \cap (y + A)| \, dx \, dy \right) E_n.$$

The result now follows as before. □

LEMMA 4. *For any interval $I = [a, b]$, $\int_{I \times I} |(x + I) \cap (y + I)| \, dx \, dy = \frac{2}{3}|I|^3$.*

PROOF. Both sides are unchanged if we replace $I$ with $I - a = [0, c]$, where $c = b - a$. So $\int_{I \times I} |(x + I) \cap (y + I)| \, dx \, dy = \int_0^c \int_0^c (c - |x - y|) \, dx \, dy = 2 \int_0^c \int_0^y (c - y + x) \, dx \, dy = \frac{2}{3} c^3$. □

THEOREM 5. *Let $A = \{x \in \mathbb{R}^2 : r(1-\varepsilon) \leq \|x\|_\infty \leq r\}$. Then $\mathcal{G}_A$ almost surely has no infinite component when $|A| \leq 1.014$.*

PROOF. Define the intervals $I = [r(1-\varepsilon), r]$ and $J = [-r, r]$, and the rectangle $S = I \times J$. Then $A \supseteq S \cup -S$. Now

$$\int_{A \times A} |(x+A) \cap (y+A)| \, dx \, dy \geq 6 \int_{S \times S} |(x+S) \cap (y+S)| \, dx \, dy,$$

where we have used symmetry to bound the integral over $x, y \in \pm S$ of $|(x + (S \cup -S)) \cap (y + (S \cup -S))|$ by six equal terms given by the six choices of sign in $x \in \pm S$, $y \in \pm S$, $(x \pm S)$ and $(y \pm S)$ that give the same contribution as $\int_{S \times S} |(x+S) \cap (y+S)| \, dx \, dy$. However, $S$ is a product $I \times J$, so we can separate variables to obtain

$$\int_{A \times A} |(x+A) \cap (y+A)| \, dx \, dy$$

$$\geq 6 \left( \int_{I \times I} |(x+I) \cap (y+I)| \, dx \, dy \right) \left( \int_{J \times J} |(x+J) \cap (y+J)| \, dx \, dy \right)$$

$$= 6(\tfrac{2}{3}|I|^3)(\tfrac{2}{3}|J|^3) \geq \tfrac{1}{24}|A|^3,$$

where we have used Lemma 4 and the bound $|A| \leq 4|I||J|$. Hence if $|A|^3 \times (1 - \tfrac{1}{24}) < 1$, then there is no percolation. This occurs if $|A| \leq 1.014$. □

**3. Upper bound.** We prove the upper bound of Theorem 1, namely that $n_c = 1 + o(\varepsilon^{1/15})$ as $\varepsilon \to 0$. The strategy is to compare the percolation process with an oriented bond percolation on $\mathbb{Z}^2$. Fix a large $R$ and integers $n < N$. Partition $\mathbb{R}^2$ into $6R \times 6R$ squares, and let the site $x \in \mathbb{Z}^2$ correspond to the square $S_x = 6Rx + [-3R, 3R]^2$ in $\mathbb{R}^2$. Our bonds $xy$ in $\mathbb{Z}^2$ correspond to certain good events in the corresponding $6R \times 12R$ rectangle $S_x \cup S_y$ of $\mathbb{R}^2$.

Roughly speaking, these good events comprise the ability to get from some set of $n$ points near the middle of $S_x$ to at least $n$ points in the middle of $S_y$ by paths that lie entirely within $S_x \cup S_y$. However, we insist on more. We "test" annuli around at most $N$ points when we construct the joining paths. Also, there are up to $3N$ points in $S_x \cup S_y$ that have been tested when constructing earlier bonds, so we wish to avoid these points and the annuli about them.

Fix sets of points $P = \{x_1, \ldots, x_n\}$ and $Q = \{y_1, \ldots, y_K\}$, $K \leq 3N$, in the rectangle $S_x \cup S_y$ that corresponds to bond $xy$ of $\mathbb{Z}^2$. We assume that no two points of $P \cup Q$ lie within $r\sqrt{\varepsilon}$ of each other. Assume further that $x_1, \ldots, x_n$ lie in the middle $4R \times 4R$ square $6Rx + [-2R, 2R]^2$ of $S_x$. Suppose $y$ is one of the two points $x + (1, 0)$ or $x + (0, 1)$. We construct subsets $P'$ and $Q'$ of the Poisson process with the following properties:

(a) $P' \subseteq 6Ry + [-2R, 2R)^2$ and $|P'| \leq n$.

(b) $Q' \subseteq S_x \cup S_y$, every element of $Q'$ is at distance at least $r$ from the boundary of $S_x \cup S_y$ and $|Q'| \leq N$.

(c) All points of $P \cup Q \cup P' \cup Q'$ are distinct and are at least $r\sqrt{\varepsilon}$ from each other.

(d) $(P' \cup Q') \cap \bigcup_{y \in Q} A_\varepsilon(y, r) = \varnothing$.

(e) Every $x' \in P'$ is joined to some point $x \in P$ by a sequence of points $y'_1, \ldots, y'_k$ of $Q'$.

We declare the oriented bond $\overrightarrow{xy}$ to be *open with respect to P and Q* if, in this construction, $|P'| = n$.

The idea is that $P$ is a set of points in $S_x$ that we can get to from the origin, and we wish to find a set of points $P'$ in $S_y$ that we can get to from $P$. The set $Q'$ consists of all the points we need to look at when constructing paths from $P$ to $P'$ and $Q$ is the set of points we have looked at previously. When constructing paths from $P$ to $P'$ it is important that we avoid annuli about points in $Q$, since this would introduce uncontrollable dependencies on our previous bonds.

Note that the openness of $\overrightarrow{xy}$ depends on the choice of $P$ and $Q$ as well as the restriction of the Poisson process to the rectangle $S_x \cup S_y$. The sets $P$ and $Q$ depend on the construction of previous bonds, which will introduce dependencies on the bonds. The key element of the proof is the following:

PROPOSITION 6. *Fix $\varepsilon > 0$ small and write $|A| = 1 + \eta$, where $\eta = \eta(r, \varepsilon)$. If $r$ is large enough that $\varepsilon \leq c\eta^{14}|\log \eta|^{-10}$, then there exist $N$, $n$ and $R$ such that the probability of a bond being open with respect to $P$ and $Q$ is at least $0.9$.*

The proof of this proposition is deferred until the next section.

PROOF OF THE UPPER BOUND OF THEOREM 1. Using Proposition 6 we complete the proof. We define the oriented percolation on $\mathbb{Z}^2$. Order the bonds in the first quadrant of $\mathbb{Z}^2$ by their $l_1$ distance from the origin and, for each distance $k$, order the bonds at distance $k$ from the origin as

$$(0,k)(0,k+1) \qquad (0,k)(1,k) \qquad (1,k-1)(1,k) \qquad \cdots \qquad (k,0)(k+1,0).$$

Suppose the bonds in this order are $\{e_1, e_2, \ldots\}$. We declare some bonds open in such a way that if there is an infinite directed open path in $\mathbb{Z}^2$ from $(0,0)$, then (with positive probability) there is an infinite path in $\mathcal{G}_A$. To this end, we inductively define, for each edge $e_i$, a subset $Q_i$ of the Poisson process with $Q_{i-1} \subseteq Q_i$ and, for each vertex $x$ joined by an open path to the origin in $\mathbb{Z}^2$, a subset $P_x$ of the Poisson process inside $S_x$.

Initially set $Q_0 = \varnothing$ and set $P_{(0,0)}$ to be any subset of $n$ points of the Poisson process that lie in $[-2R, 2R)^2$ and such that no two points of $P_{(0,0)}$ lie within $r\sqrt{\varepsilon}$ of each other. It is easy to check that such a set exists with high probability.

Now suppose we have defined the openness of the bonds $e_j$ for $j < i$ and the set $Q_i$ and all relevant $P_x$. We now consider the bond $e_i = xy$.

(a) If there is no directed path consisting of open bonds from $(0,0)$ to $x$, set $xy$ to be open and $Q_{i+1} = Q_i$.

(b) If there is a directed open path from $(0,0)$ to $x$, declare $xy$ to be open if it is open with respect to $P = P_x$ and $Q = Q_i \cap (S_x \cup S_y)$. Set $Q_{i+1} = Q_i \cup Q'$ and, if $xy$ is open and $P_y$ is not yet defined, set $P_y = P'$.

Condition (a) is a technical condition which clearly does not affect whether or not $(0,0)$ is in an infinite cluster.

By (b), at each step $|Q_i \cap S_z| \leq kN$, where $k$ is the number of edges $e_j$, $j < i$, meeting $z$. Thus, given the ordering of bonds as above, if $e_i$ is a vertical bond, $|Q_i \cap S_x| \leq 2N$ and $|Q_i \cap S_y| \leq N$, whereas if $e_i$ is a horizontal bond, $|Q_i \cap S_x| \leq 3N$ and $|Q_i \cap S_y| = 0$. Hence the set $Q$ in (b) always satisfies $|Q| \leq 3N$.

There are two edges $e_i = xy$ with a given value of $y$, so there are two chances for $P_y$ to be defined in (b). Clearly $P_y$ is defined iff $y$ is joined to $(0,0)$ by an open path. Also, by construction, all points in $(\bigcup_i Q_i) \cup (\bigcup_x P_x)$ are at least $r\sqrt{\varepsilon}$ from each other. If $x$ is joined to $(0,0)$ by an open path $x = x_k, \ldots, x_0 = (0,0)$, then each point in $P_{x_{i+1}}$ is joined by a path in $\mathcal{G}_A$ to a point in $P_{x_i}$. Hence there is a path from any point of $P_x$ to one of the $n$ points of $P_{(0,0)}$.

Finally, by Proposition 6, if $\varepsilon \leq c\eta^{14} |\log \eta|^{-10}$ [which is satisfied if $\eta = \Omega(\varepsilon^{1/15})$], then each bond is open with probability bounded below by 0.9 even conditioned on the state of all previous edges and regions of $\mathbb{R}^2$ on which they depend (the annuli around the points of $Q_i$). Thus the percolation process on $\mathbb{Z}^2$ stochastically dominates an independent oriented bond percolation with bond probability 0.9. However, $(0,0)$ is then in an infinite cluster with positive probability [see, e.g., Balister, Bollobás and Stacey (1994)]. The result now follows. □

**4. Proof of Proposition 6.** The proof of the bound is complicated by the fact that the annuli intersect, so we first consider the simpler case when we ignore these intersections and model the percolation by a branching process [see, e.g., Athreya and Ney (1972)]. We generally refer to the points of a branching process as *nodes* to avoid confusion with the points of our Poisson process.

LEMMA 7. *Consider a branching process where at each step each node branches into several new nodes according to independent identical Poisson distributions with mean $1+\eta$. Let $N_t$ be the number of nodes at time $t > 0$. Then $\mathbb{P}(N_t \geq (1+\eta)^t) \geq \eta \mathbb{P}(N_t = 0)^2 \geq \eta e^{-2(1+\eta)}$.*

PROOF. It is easy to show by induction on $t$ that

$$\mathbb{E} N_t = (1+\eta)^t$$

and

$$\operatorname{Var} N_t = (1+\eta)^t((1+\eta)^t - 1)/\eta < (\mathbb{E} N_t)^2/\eta.$$

Let $X_t = N_t/\mathbb{E} N_t$, so that $\mathbb{E} X_t = 1$ and $\operatorname{Var} X_t < 1/\eta$. By Cauchy–Schwarz,

$$\mathbb{E}(I_{X_t \geq 1})\mathbb{E}((X_t - 1)^2) \geq (\mathbb{E}((X_t - 1)I_{X_t \geq 1}))^2.$$

Now

$$\mathbb{E}((X_t - 1)^2) = \operatorname{Var} X_t < 1/\eta,$$

$$\mathbb{E}((X_t - 1)I_{X_t \geq 1}) = \mathbb{E}((1 - X_t)I_{X_t \leq 1})$$

$$\geq \mathbb{P}(X_t = 0) = \mathbb{P}(N_t = 0).$$

Hence $\mathbb{P}(N_t \geq (1+\eta)^t) = \mathbb{E}(I_{X_t \geq 1}) \geq \eta \mathbb{P}(N_t = 0)^2$. Finally $\mathbb{P}(N_t = 0) \geq \mathbb{P}(N_1 = 0) = e^{-(1+\eta)}$. □

LEMMA 8. *Assume $\eta > 0$ is sufficiently small and consider a branching process where at each step each node branches into several new nodes independently according to identical Poisson distributions with mean $1 + \eta$. Suppose also that we remove nodes so there are at most $K$ nodes at each step and assume $T \leq e^{\eta K} \eta/3$. Then the probability that there is at least one node at time $T$ is at least $\eta/3$.*

PROOF. Let $N_t$ be the number of nodes at time $t$ and consider the random variable $X_t = \exp(-\lambda N_t)$, where $\lambda$ is the positive solution to the equation $(1 - e^{-\lambda})(1 + \eta) = \lambda$ [unique since $(1 - e^{-\lambda})/\lambda$ monotonically decreases from 1 to 0]. Now

$$N_{t+1} = \min(N'_{t+1}, K) \qquad \text{where } N'_{t+1} = \sum_{i=1}^{N_t} Y_i$$

and $Y_i$ are independent identical Poisson random variables of mean $1 + \eta$. Also

$$\mathbb{E}\exp(-\lambda Y_i) = e^{-(1+\eta)} \sum_{n=0}^{\infty} \frac{(1+\eta)^n}{n!} e^{-\lambda n} = \exp((e^{-\lambda} - 1)(1+\eta)) = e^{-\lambda}.$$

Hence
$$\mathbb{E}(\exp(-\lambda N'_{t+1})|N_t) = \prod_{i=1}^{N_t} \mathbb{E}\exp(-\lambda Y_i) = \exp(-\lambda N_t) = X_t,$$
but if $N_{t+1} \neq N'_{t+1}$, then $K = N_{t+1} < N'_{t+1}$. Hence
$$0 \leq \mathbb{E}(\exp(-\lambda N_{t+1}) - \exp(-\lambda N'_{t+1})|N_t) \leq e^{-\lambda K}.$$
So $\mathbb{E}(X_{t+1} \mid N_t) \leq X_t + e^{-\lambda K}$ and thus $\mathbb{E}X_{t+1} \leq \mathbb{E}X_t + e^{-\lambda K}$. However, $\mathbb{E}X_t \geq \mathbb{P}(N_t = 0)$ and $\mathbb{E}X_0 = \exp(-\lambda N_0) = e^{-\lambda}$. Hence
$$\mathbb{P}(N_T = 0) \leq \mathbb{E}X_T \leq e^{-\lambda} + Te^{-\lambda K}.$$
It remains to bound $\lambda$. Let $f(x) = (1 - e^{-x})(1 + \eta) - x$, so that $\lambda$ is a solution of $f(\lambda) = 0$. Then $f(\eta) = (\eta - \eta^2/2 + O(\eta^3))(1+\eta) - \eta$ is positive for small $\eta$, but $f(1 + \eta) < 0$. Hence $\lambda > \eta$ and thus $\mathbb{P}(N_t = 0) \leq e^{-\eta} + Te^{-\eta K} \leq e^{-\eta} + \eta/3$, which is at most $1 - \eta/3$ for sufficiently small $\eta$. □

We now consider a simplified version of our percolation process in which each step is independent of all previous steps. Consider a branching process where at each step each node branches into several new nodes according to independent identical Poisson distributions with mean $1 + \eta$. Assign to each node $v$ (other than the root node) a random variable $A_v$ uniformly distributed in $A = A_\varepsilon(0, r)$ independently of the branching process and all other $A_u$'s. Fix a position $z_0$ in $\mathbb{R}^2$ for the root node and define the position $z_v$ of a node $v$ to be $z_0 + \sum_u A_u$, where the sum runs over all ancestors $u$ of $v$ back to the root node. Let $\mathcal{T}_A$ be the random graph with vertices $z_v$ and edges $z_v z_{v'}$ for all child nodes $v'$ of $v$. Also, define $\mathcal{T}_A^t$ to be the set of nodes that occur at time $t$. The process $\mathcal{T}_A$ clearly approximates the percolation process $\mathcal{G}_A$, but it differs in that the distribution of points (child nodes) in $A_\varepsilon(z_v, r)$ is independent of the process up to that point, whereas in $\mathcal{G}_A$ the points in $A_\varepsilon(z_v, r)$ depend on points in previous annuli where they intersect.

LEMMA 9. *Pick $z_0 \in [-2R, 2R)^2$ and consider the branching process $\mathcal{T}_A$ defined above, except that at each step (if necessary) we remove nodes (randomly, independent of their locations) so there are at most $K$ new nodes at each step. Run this process up to time $T = (R/r)^2$. Define the event $\mathcal{E}$ to be the event that the process has not died out at time $T$ (i.e., $\mathcal{T}_A^T \neq \varnothing$) and that, picking a node at random from $\mathcal{T}_A^T$, this node lies in the square $(0, 6R) + [-R, R)^2$ and all its ancestors lie in $[-3R + r, 3R - r) \times [-3R + r, 9R - r)$. Then there exists an absolute constant $c_0 > 0$, independent of $\eta$, $\varepsilon$, $r$ and $R$, such that for $\eta$ sufficiently small and $R/r$ sufficiently large, $\mathbb{P}(\mathcal{E}|\mathcal{T}_A^T \neq \varnothing) \geq c_0$.*

PROOF. Conditioning on $\mathcal{T}_A^T \neq \varnothing$, pick a node $v$ at random from $\mathcal{T}_A^T$. Note that the choice of $v$ is independent of the locations of the nodes. Consider the locations of the nodes on the path from the root node to $v$. These form a random walk with steps taken uniformly from $A = A_\varepsilon(0, r)$. The probability we are interested in is bounded below by the probability of this random walk meeting $(0, 6R) + [-R, R]^2$ before hitting the boundary of the rectangle $[-3R + r, 3R - r) \times [-3R + r, 9R - r)$ and before time $T$. Each step has mean 0 and variance in either coordinate direction of $cr^2$, where $c = c(\varepsilon)$ can be bounded above and below by positive constants independently of $\varepsilon$. By scaling the dimensions by $1/R$ and time by $c'(r/R)^2$, this random walk for large $R/r$ can be approximated by a two-dimensional Brownian motion, run for constant time $t_0$, starting in $[-2, 2]^2$. The probability of the Brownian motion lying in $[-0.9, 0.9] \times [5.1, 6.9]$ at time $t_0$ without hitting the boundary of $[-2.9, 2.9] \times [-2.9, 8.9]$ before time $t_0$ is bounded below by a constant $c'' > 0$. Choosing $c_0 < c''$ we see that provided $R/r$ is larger than some absolute constant, the probability that $v$ lies in the square $(0, 6R) + [-R, R]^2$ and all its ancestors lie in the rectangle $[-3R + r, 3R - r) \times [-3R + r, 9R - r)$ is at least $c_0$. Moreover, the bound on $R/r$ and the value of $c_0$ can be chosen independently of $r$, $R$, $\eta$ and $\varepsilon$. □

We now need to deal with the dependencies caused by the intersections of annuli. We start with a couple of geometric lemmas.

LEMMA 10. *There exists a constant $c_1 > 0$ such that for any $r$, $\varepsilon$ and $x$, $y$, with $\|x - y\|_2 \geq r\sqrt{\varepsilon}$, we have $|A_\varepsilon(x, r) \cap A_\varepsilon(y, r)| \leq c_1 |A| \sqrt{\varepsilon}$.*

PROOF. Tedious, but straightforward verification. Note that there are two cases when the bound is tight: one when $\|x - y\|_2 \approx r\sqrt{\varepsilon}$ and the other when $\|x - y\|_2 \approx r(2 - \varepsilon)$. □

Note that this lemma fails for the square annulus when $\|x - y\|_\infty \approx r(2 - \varepsilon)$.

LEMMA 11. *Let $x_1, \ldots, x_k$ be points, no two of which are within $r\sqrt{\varepsilon}$ of each other. Let $A_i = A_\varepsilon(x_i, r)$ and let $B_i = A_1(x_i, r\sqrt{\varepsilon})$ be the ball around $x_i$ of radius $r\sqrt{\varepsilon}$. Then*

$$\left| A_i \cap \left( \bigcup_{j \neq i} (A_j \cup B_j) \right) \right| \leq c_2 k |A| \sqrt{\varepsilon}.$$

PROOF. The region $A_i \cap B_j$ has area $O((r\varepsilon)(r\sqrt{\varepsilon})) = O(|A|\sqrt{\varepsilon})$ since $A_i$ is of "width" $r\varepsilon$ and $B_j$ is of diameter $2r\sqrt{\varepsilon}$. The region $A_i \cap A_j$ has area $O(|A|\sqrt{\varepsilon})$ by Lemma 10. The result follows. □

PROOF OF PROPOSITION 6. The strategy of the proof is to find some points $x_i''$ in $6Ry + [-R, R]^2$ by comparing the process $\mathcal{G}_A$ with $\mathcal{T}_A$. There will, however, be rather fewer than $n$ such points, so we then run the percolation for $R/r$ further steps in $\mathcal{G}_A$ to obtain sufficiently many points in $6Ry + [-2R, 2R]^2$ (note that we cannot travel more than distance $R$ in $R/r$ steps).

Pick each $x_i \in P$ in turn and run a truncated branching process $\mathcal{T}_A$ as in Lemma 9, starting at $x_i$. We call a node *bad* if it lies in any annulus $A_\varepsilon(z, r)$ or ball $A_1(z, r\sqrt{\varepsilon})$, where $z$ is in $Q$, $P \setminus \{x_i\}$ or any one of the nodes of the branching process other than the parent of $v$, or any one of the points in the branching processes already constructed for $x_j \in P$, $j < i$. There are at most $3N + n + nK(R/r)^2$ such values of $z$. Set

$$K = \frac{1}{\eta} \log\left(\frac{3R^2}{\eta r^2}\right). \tag{1}$$

It is clear that we can run the same branching process in the percolation model $\mathcal{G}_A$, coupling $\mathcal{G}_A$ with $\mathcal{T}_A$ so that they agree up until we hit a bad node of $\mathcal{T}_A$. By Lemma 8, the branching process has not died out by time $T = (R/r)^2 = e^{\eta K} \eta/3$ with probability at least $\eta/3$. Moreover, the probability that a given node in this process is bad conditioned on all its predecessors being good is at most $c_2(3N + n + nK(R/r)^2)|A|\sqrt{\varepsilon}/|A|$ by Lemma 11. Also this is independent of any event involving the existence of its decendants in the branching process. Thus, if we condition on $\mathcal{T}_A^T \neq \varnothing$ and we pick a path of length $T$ (independently of locations of the nodes), the probability that this path contains a bad node is at most $(R/r)^2 c_2(3N + n + nK(R/r)^2)\sqrt{\varepsilon}$. We require

$$(R/r)^2 c_2(3N + n + nK(R/r)^2)\sqrt{\varepsilon} \leq c_0/2. \tag{2}$$

Thus by Lemma 9, conditioning on $\mathcal{T}_A^T \neq \varnothing$, we obtain a point $x_i''$ in $6Ry + [-R, R]^2$ joined to $x_i$ in the percolation process within $S_x \cup S_y$ with probability at least $c_0 - c_0/2 = c_0/2$. Thus, by Lemma 8, the probability of finding such a point $x_i''$ is at least $(c_0/2)(\eta/3) = c_0\eta/6$, independently of all previously found points. The number of points $x_i''$ found is hence stochastically bounded below by a binomial variable $X$ with mean $c_0 n\eta/6$. Now take the points $V_0 = \{x_1'', \ldots, x_X''\}$ and construct inductively sets $V_i$ by taking all points joined to some point in $V_{i-1}$ which are not in $A_\varepsilon(z, r) \cup A_1(z, r\sqrt{\varepsilon})$ for any point $z$ already considered. Repeat for $R/r$ steps or until $|V_i| \geq n$ if this occurs earlier. All points of $V_i$ must then lie in $6Ry + [-2R, 2R]^2$ and the cardinalities $|V_i|$ are stochastically dominated by a Poisson branching process with mean $1 + \eta/2$ provided

$$c_2(3N + n + nK(R/r)^2 + n(R/r))\sqrt{\varepsilon} \leq \eta/2. \tag{3}$$

CONTINUUM PERCOLATION 11

Define $Q'$ to be all the new points encountered above which are in $6Ry + [-3R, 3R]^2$ except those of the last $V_i$ and let $P'$ consist of $n$ points from this last $V_i$, if they exist. Let

$$N = nK(R/r)^2 + n(R/r). \tag{4}$$

Then $N > |Q'|$. Finally, the probability of the bond being open is bounded below by the probability that a branching process with parameter $1 + \eta/2$ starting from $X$ points will have at least $n$ points by time $R/r$. If

$$n \leq (1 + \eta/2)^{R/r}, \tag{5}$$

then by Lemma 7, the probability that such a branching process starting with one node has at least $n$ nodes by time $R/r$ is at least $(\eta/2)e^{-2(1+\eta/2)}$, which is at least $\eta/20$ for small $\eta$. Now the probability of not having at least $n$ nodes by time $R/r$ starting from $X$ points is bounded above by $(1 - \eta/20)^X$, which has expected value

$$\sum_i \binom{n}{i} (c_0\eta/6)^i (1 - c_0\eta/6)^{n-i} (1 - \eta/20)^i$$
$$= (1 - c_0\eta^2/120)^n$$
$$\leq \exp(-c_0 n\eta^2/120).$$

The result follows provided

$$c_0 n\eta^2/120 \geq -\log(0.1). \tag{6}$$

It remains to choose the parameters $\varepsilon$, $n$, $N$ and $R$ so that equations (1)–(6) are satisfied. For (6) we can define $n = c_3 \eta^{-2}$ for some $c_3 > 0$. Then (5) is satisfied if we take $R/r = c_4 \eta^{-1} |\log \eta|$. From (1) and (4) we get $N = O(\eta^{-5}|\log \eta|^3)$, so (2) and (3) are now satisfied if $\sqrt{\varepsilon} \leq c_5 \eta^7 |\log \eta|^{-5}$ or, equivalently, $\varepsilon \leq c\eta^{14} |\log \eta|^{-10}$. $\square$

Note that the sets $P'$, $Q'$ depend only on the points of the Poisson process inside the region $(S_x \cup S_y) \cap (\bigcup_{z \in P \cup Q'} A_\varepsilon(z,r)) \setminus (\bigcup_{z \in Q} A_\varepsilon(z,r))$, so each bond is independent of the regions tested when constructing previous bonds.

DEPARTMENT OF MATHEMATICS
UNIVERSITY OF MEMPHIS
DUNN HALL
3725 NORISWOOD
MEMPHIS, TENNESSEE 38152
USA
E-MAIL: balistep@msci.memphis.edu
E-MAIL: bollobas@msci.memphis.edu
E-MAIL: mjw1009@cam.ac.uk